\documentclass{article}

\usepackage{amsmath}

\usepackage{amsthm}

\begin{document}

\begin{center}
\large \textbf{p-Parts of Stabilizers in Primitive Permutation Groups}
\end{center}

\bigskip
\bigskip

\begin{center}
\textbf{ABSTRACT}
\end{center}

\medskip

\noindent Let $(G,\Omega)$ be a finite primitive permutation group.
Let $p^2$ divide $|G|$, for a prime $p$. We show that when $G$ is
solvable, there exists a subset of $\Omega$ whose stabilizer $S$
has the property that $1<|S|_p<|G|_p.$ We offer a counting argument
which should be helpful when $G$ is not solvable.

\bigskip
\bigskip
\bigskip
\bigskip
\bigskip

\begin{center}
David Gluck\\
Department of Mathematics\\
Wayne State University\\
Detroit, Michigan 48202\\
d.gluck@me.com\\
\end{center}

\bigskip
\bigskip
\bigskip
\bigskip
\bigskip
\bigskip
\bigskip
\bigskip
\bigskip
\bigskip
\bigskip
\bigskip
\bigskip
\bigskip

\noindent Key words: primitive permutation group, affine group, set-stabilizer

\bigskip
\bigskip
\bigskip
\bigskip

\noindent MSC: 20B15, 20C20

\newpage

\begin{center}

\textbf{0. Introduction}

\end{center}

A key step in the recently completed proof [8] of Brauer's
height zero conjecture is the so-called Gluck-Wolf-Navarro-Tiep
theorem [6,10], which asserts that if $N$ is a normal subgroup
of a finite group $G$, $\theta \in \mathrm{Irr}(N)$, and 
$(p, \chi(1)/\theta(1))=1$ for a prime $p$ and all $\chi \in \mathrm{Irr}
(G|\theta)$, then $G/N$ has abelian Sylow $p$-subgroups. The
proof of this theorem involved the analysis of certain imprimitive 
modules, and therefore required the classification of certain
permutation groups, which was done in [4, Theorem 2].

More precisely, the authors of [4] defined a permutation group
$(G, \Omega)$ to be $p$-concealed if $p$ divides $|G|$ and every
subset of $\Omega$ is stabilized by some Sylow $p$-subgroup of $G$.
Apart from some examples in which $|\Omega|\leq 8$, the primitive 
$p$-concealed groups are certain alternating and symmetric groups
with natural permutation domain $\Omega$.

The Eaton-Moret\'{o} conjecture [3], a vast generalization of Brauer's
height zero conjecture, leads to analogous, but more difficult classification
problems for permutation groups. One case of the Eaton-Moret\'{o}
conjecture gives rise to the problem of classifying the primitive permutation
groups $(G,\Omega)$ with the property that for each subset $\Delta$ of
$\Omega$, the setwise stabilizer $\mathrm{Stab}_G(\Delta)$ has either
$p^{\prime}$-order or $p^{\prime}$-index in $G$. Informally, we call
such permutation groups $p$-extreme. If $p^2$ does not divide $|G|$,
we say that $(G, \Omega)$ is trivially $p$-extreme. Clearly $p$-concealed groups
are $p$-extreme, and so certain alternating and symmetric groups are
nontrivial primitive $p$-extreme groups; the smallest of these is $\mathrm{Alt}(7)$
with $p=2$.

These are possibly the only nontrivial primitive $p$-extreme groups.
Our main result, Theorem A, says that there are no nontrivial solvable
primitive $p$-extreme groups. It seems likely that there are also no
nontrivial $p$-solvable primitive $p$-extreme groups. 

Actually, it is easier to work with the complementary notion of
$p$-moderate (that is, not $p$-extreme) permutation groups. Our Theorem A
then states formally that every primitive solvable permutation group
$(G,\Omega)$, with $|G|$ divisible by $p^2$, is $p$-moderate. A
regular orbit theorem of Dolfi [2, Theorem 3.4] plays a major role in
the proof of Theorem A for odd $p$. When $p=2$, the argument using
Dolfi's theorem breaks down so instead we use regular orbit results of
Yang and his collaborators [7, 12]. In a brief final section, we offer an
elementary counting argument, Proposition 3.1, which should be helpful
in showing that most non-solvable primitive permutation groups are
$p$-moderate.

While our motivation comes from representation theory, there are other reasons
to be interested in stabilizers in finite permutation groups. The study of such 
stabilizers is still in its early stages, but there has been significant recent
activity on this topic; see, for example [1, 5, 11].

I thank Juan Martinez, Alex Moret\'{o}, and Gabriel Souza for helpful
communications related to this paper. I thank the referee for a careful reading 
of this paper, which has led to improvements in the presentation.

\bigskip

\begin{center}

\textbf{1. A reduction theorem}

\end{center}

\bigskip

\noindent \textit{Notation and definitions.} An affine permutation group
$(G, \Omega)$ is one in which $\Omega$ is a finite vector space $V$
and $G$ is a semidirect product $VH$, where $H$ acts faithfully on $V$
by linear transformations and $V$ acts regularly on itself by translations.
Such an affine group is a primitive permutation group if and only if $V$ is
an irreducible $H$-module. We recall that every solvable primitive 
permutation group is an affine group. Following [4], we say that a permutation
group $(G, \Omega)$ is $p$-concealed for a prime $p$ if $p$ divides $|G|$ and every subset of
$\Omega$ is stabilized by some Sylow $p$-subgroup of $G$. 

\medskip

\noindent \textbf{Proposition 1.1} A $p$-solvable primitive permutation group
$(G, \Omega)$ is $p$-concealed if and only if $(G, \Omega)$ is an affine
group and $(p, G, |\Omega|)$ is $(2, D_6, 3), (2,D_{10}, 5)$, or $(3, J, 8)$,
where $J$ is the full affine semilinear group over $GF(8)$; $J$ is solvable of
order 168.

\medskip

\noindent \textit{Proof.} This follows from [4, Theorem 2]. \qedsymbol

\medskip

\noindent \textit{Definition}. We define a permutation group $(G,\Omega)$ to
be $p$-moderate if there exists a subset $\Delta$ of $\Omega$ whose setwise
stabilizer $\mathrm{Stab}_G(\Delta)$satisfies
$$1<|\mathrm{Stab}_G(\Delta)|_p<|G|_p.$$

\medskip

\noindent \textbf{Proposition 1.2} Let $G=VH$ be a (not necessarily primitive)
$p$-solvable affine permutation group for a prime $p$. Suppose that
$(p, |V|)=1$ and $p^2$ divides $|H|$. Let $V=V_1\oplus \ldots \oplus V_k$,
where each $V_i$ is an irreducible $H$-module. Suppose for each $i$ that
the primitive affine group $G_i:=V_i(H/C_H(V_i))$ is either $p$-moderate
or that $p^2$ does not divide $|G_i|$. Then the affine group $(G,V)$ is 
$p$-moderate.

\medskip

\noindent \textit{Proof.} We proceed by induction on $k$. The result is
clear if $k=1,$ so we assume that $k>1.$ Suppose first that $p^2$ divides,
say, $|G_1|.$ Let $\Delta_1 \subseteq V_1$ with
$$1<|\mathrm{Stab}_{G_1}(\Delta_1)|_p<|G_1|_p.$$
Since $\mathrm{Stab}_G(\Delta_1)=\mathrm{Stab}_{VH}(\Delta_1)$ and
$p$ does not divide $|V_1|$, we may, after replacing $\Delta_1$ by a
$V_1H$-conjugate subset of $V_1$, assume that a Sylow $p$-subgroup of 
$\mathrm{Stab}_G(\Delta_1)$ is contained in $H$. Then $|\mathrm{Stab}_G(\Delta_1)|_p
=|\mathrm{Stab}_H(\Delta_1)|_p.$ Viewing $\Delta_1$ as a subset of $V$,
we have
$$|\mathrm{Stab}_G(\Delta_1)|_p=|C_H(V_1)|_p|\mathrm{Stab}_{H/C_H(V_1)}(\Delta_1)|_p
=|C_H(V_1)|_p|\mathrm{Stab}_{G_1}(\Delta_1)|_p.$$ Hence
$$1<|\mathrm{Stab}_G(\Delta_1)|_p<|C_H(V_1)|_p|G_1|_p=|G|_p$$ and so $(G,V)$ is 
$p$-moderate, as desired.

Hence $|G_i|_p\leq p$ for $1\leq i \leq k$. Suppose that $|G_i|_p=p$ for $1 \leq i
\leq m$ and $|G_i|_p=1$ for $m+1 \leq i \leq k$. Since $p^2$ divides $|G|$, we have 
$m \geq 2.$ We claim that $|H/C_H(V_1 \oplus V_j)|_p=p^2$ for some $j$ with 
$2 \leq j \leq m$. Suppose the contrary. Then
$$|H/C_H(V_1 \oplus V_j)|_p=|H/C_H(V_1)|_p=p$$
for $2 \leq j \leq m$. It follows that $|C_H(V_1)|_p=|C_H(V_1 \oplus V_j)|_p$ for
$2 \leq j \leq m$. Let $Q \in \mathrm{Syl}_p(C_H(V_1))$. Since $C_H(V_1 \oplus V_j)$
is a normal subgroup of $C_H(V_1)$, we have
$$Q \cap C_H(V_1 \oplus V_j) \in  \mathrm{Syl}_p(C_H(V_1 \oplus V_j)),$$
and so
$$|Q \cap C_H(V_1 \oplus V_j)|=|H|_p/p=|C_H(V_1)|_p=|Q|.$$
Thus $Q \leq C_H(V_1 \oplus V_j).$ Hence $Q$ centralizes $V_j$ for 
$2 \leq j \leq m.$ By definition of $m, Q$ also centralizes $V_j$ for  $m+1 \leq j
\leq k.$ Hence $Q$ centralizes $V$, against the fact that $H$ acts faithfully on
$V$. This proves the claim.

After renumbering, we may assume that $|H/C_H(V_1 \oplus V_2)|_p=p^2.$
If $k>2,$ the inductive hypothesis, applied to the affine group $G_{12}:=
(V_1 \oplus V_2)(H/C_H(V_1 \oplus V_2)),$ yields a subset $\Delta$ of 
$V_1 \oplus V_2$ such that $$|\mathrm{Stab}_{G_{12}}(\Delta)|_p=p<p^2
=|G_{12}|_p.$$ Since $\mathrm{Stab}_G(\Delta) \leq (V_1\oplus V_2)H,$
it follows that
$$|\mathrm{Stab}_G(\Delta)|_p=p|C_H(V_1 \oplus V_2)|_p<|H|_p=|G|_p.$$
Hence the affine group $G$ is $p$-moderate.

Thus we may assume that $k=2, V=V_1 \oplus V_2, |H/C_H(V_i)|_p=p$
for $i=1,2,$ and $|H|_p=|G|_p=p^2.$ For $G_1$ and $G_2$ as above,
suppose that $(G_1, V_1)$, say, is not $p$-concealed. Choose $\Delta _1
\subseteq V_1$ such that $|\mathrm{Stab}_{G_1}(\Delta _1)|_p=1.$ Viewing $\Delta_1$
as a subset of $V$, we have $C_H(V_1) \leq \mathrm{Stab}_G(\Delta_1) \leq
V_1H,$ and so $|\mathrm{Stab}_G(\Delta_1)|_p=|C_H(V_1)|_p=p.$ Hence 
$(G,V)$ is $p$-moderate.

Thus we may assume that both $(G_1, V_1)$ and $(G_2, V_2)$ are 
$p$-concealed. By Proposition 1.1, $p$ is 2 or 3. If $p=2,$ then either 
$|V_1|=|V_2|=3$ and $G_1$ and $G_2$ are dihedral of order 6, or
$|V_1|=|V_2|=5$ and $G_1$ and $G_2$ are dihedral of order 10.
If $p=3,$ then $|V_1|=|V_2|=8$ and $G_1$ and $G_2$ are isomorphic
to $J$, the full affine semilinear group over $GF(8).$

If $p=2$ and $|V_1|=|V_2|=3,$ then, since $|G|_2=4, G$ must be
the full direct product $D_6 \times D_6.$ Let
$$\Delta =\{(0,0),(1,1)\} \subseteq V_1 \oplus V_2=V.$$
Now $G$ consists of all transformations $T(a,b)D(\epsilon _1,
\epsilon_2)$ where $(a,b) \in V_1 \oplus V_2$ and $\epsilon_i 
\in \{\pm 1\}$ for $i=1,2.$ Here $T(a,b)$ is the affine translation
by the vector $(a,b)$ and $D(\epsilon _1, \epsilon _2)$ is the diagonal linear
transformation on $V$ which takes $(v_1, v_2)$ to $(\epsilon_1v_1,
\epsilon_2v_2).$ Suppose that $T(a,b)D(\epsilon_1,\epsilon_2)$ stabilizes
$\Delta$. Then either $(\epsilon_1a, \epsilon_2b)=(1,1)$ and
$(\epsilon_1(a+1), \epsilon_2(b+1))=(0,0)$ or $(\epsilon_1a, \epsilon_2b)=
(0,0)$ and $(\epsilon_1(a+1), \epsilon_2(b+1))=(1,1).$ The first possibility
yields $a=b=-1$ and $\epsilon_1=\epsilon_2=-1.$ The second possibility yields 
$a=b=0$ and $\epsilon_1=\epsilon_2=1.$ Hence $|\mathrm{Stab}_G(\Delta)|
=2,$ and so $G$ is $2$-moderate, as desired. The same construction
works when $p=2$ and $|V_1|=|V_2|=5.$

Suppose next that $p=3, |V_1| =|V_2|=8,$ and $G_1$ and $G_2$ are 
isomorphic to $J$, the affine semilinear group over $GF(8).$ Then $G$
has index 1 or 7 in $J\times J$, an affine group on $GF(8) \oplus GF(8).$
We note that $J$ has 28 Sylow  3-subgroups, each of which has 4 orbits on
$GF(8)$, of sizes 1,1,3, and 3. Hence $J \times J$ has $28^2$ Sylow 3-subgroups,
each of which has 16 orbits on $GF(8) \oplus GF(8).$ It follows that a
Sylow 3-subgroup of $G$ stabilizes $2^{16}$ subsets of $V$, and so at most 
$(28^2)(2^{16})$ subsets of $V$ are stabilized by some Sylow 3-subgroup of
$G$. Now let $t\in G$ be an element of order 3 which acts trivially on $V_2.$
Then $t$ has (4)(8)=32 orbits on $V$, and so $t$  stabilizes $2^{32}$ subsets 
of $V$. Since $(28^2)(2^{16})<2^{32},$ some subset of $V$ is stabilized by
$t$, but by no Sylow 3-subgroup of $G$. Hence $G$ is 3-moderate, as desired. \qedsymbol

\medskip

\noindent \textit{Remark.} The idea in the last paragraph of the preceding
proof is applicable much more generally; see Proposition 3.1 below.

\medskip

\noindent \textbf{Theorem A} Let $(G,\Omega)$ be a primitive solvable
permutation group. Suppose that $p^2$ divides $|G|,$ for a prime $p$.
Then $(G, \Omega)$ is $p$-moderate.

\medskip

Theorem A will be proved in the next section. The following reduction theorem
is the main result of this section. Recall that for any finite group $H,  O^{p^{\prime}}(H)$
is the smallest normal subgroup $K$ of $H$ such that $H/K$ is a $p^{\prime}$-
group, while $O_{p^{\prime},p}(H)$ is defined by 
$$O_{p^{\prime}, p}(H)/O_{p^{\prime}}(H)=O_p(H/O_{p^\prime}(H)).$$

\medskip

\noindent  \textbf{Theorem 1.3} Let $(G,\Omega)$ be a counterexample to
Theorem A with $|G|$ minimal. Let $\Omega=V$ and $G=VH$ as above. Then
\begin{enumerate}
\item $(p,|V|)=1$.
\item $H=O^{p^{\prime}}(H)$.
\item $H=O_{p^{\prime},p}(H)$.
\item $|H|_p=p^2$ and $H$ has elementary abelian Sylow $p$-subgroups.
\item $V$ is a primitive $H$-module.
\end{enumerate}

\medskip

\noindent \textit{Proof.} Suppose assertion (1) is false. If $|V|=p,$ then $H$
is isomorphic to a subgroup of $\mathrm{Aut} (V),$ which has order 
$p-1.$ Hence $p^2$ does not divide $|G|$, a contradiction. Thus $|V|_p>p.$
Let $W\leq V$ be a subgroup of order $p$. Let $\Delta=W.$ The stabilizer of
$\Delta$ in the group of affine translations of $V$ consists of translations by
vectors in $W$. Hence this stabilizer has order $p$. Let
$S \in \mathrm {Syl}_p(\mathrm{Stab}_G(\Delta)).$ Since $S$ does not
contain the full group of affine translations of $V$, we have
$S \notin \mathrm{Syl}_p(G).$ Since $|S| \geq p,$ it follows that $(G,V)$ is
$p$-moderate, a contradiction. This proves (1).

Suppose assertion (2) is false. Let $K=O^{p^{\prime}}(H).$ Then $V_K$ is
completely reducible. Say $V_K=V_1 \oplus \ldots \oplus V_k,$ with each
$V_i$ an irreducible $K$-module. For $1 \leq i \leq k,$ let $G_i
=V_i(K/C_K(V_i)),$ as in Proposition 1.2. By the minimality of $(G,V)$
as a counterexample to Theorem A, each $(G_i, V_i)$ is $p$-moderate
whenever $p^2$ divides $|G_i|$. By Proposition 1.2, $(VK, V)$ is
$p$-moderate. Hence there exists $\Delta \subseteq V$ such that
$|\mathrm{Stab}_{VK}(\Delta)|_p >1$ and no Sylow $p$-subgroup  of
$VK$ stabilizes $\Delta.$ Since $p$ does not divide $|G/K|,$ this means
that no Sylow $p$-subgroup of $G$ stabilizes $\Delta$. Hence $(G, V)$ is
$p$-moderate, a contradiction. This proves (2).

Now $H=O^{p^{\prime}}(H)$ must contain a normal subgroup $M$ with
$|H:M|=p.$ Suppose that $p^2$ divides $|M|.$ Since $V_M$ is completely reducible, the argument in the 
preceding paragraph shows that there exists a subset $\Delta$ of $V$
such that $1<|\mathrm{Stab}_{VM}(\Delta)|_p<|M|_p.$ Hence
$p \leq |\mathrm{Stab}_G(\Delta)|_p <|G|_p,$ and so $(G,V)$ is
$p$-moderate, a contradiction. We conclude that $|M|_p=p$ and so
$|H|_p=p^2.$ Since $H$ is $p$-solvable with abelian Sylow $p$-subgroups,
we have $H=O_{p^{\prime}, p, p^{\prime}}(H).$ Thus (2) implies (3).

Suppose, to get a contradiction, that $G$ has cyclic Sylow $p$-subgroups.
Choose $t \in H$ of order $p$. Choose $v \in V$such that $vt \ne v.$ Let
$\Delta=\{v, vt, \ldots, vt^{p-1}\}.$ Then $t$ stabilizes $\Delta$. Since $G$
is not $p$-moderate, $t$ lies in a Sylow $p$-subgroup $P$ of $\mathrm{Stab}_G
(\Delta)$ with $|P|=p^2=|G|_p$. Thus $P$ is also a Sylow $p$-subgroup of $G$.
Let $P=\langle u \rangle,$ with $u^p=t.$ Since $|\mathrm{Sym}(p)|_p=p,$
it follows that $u^p$ lies in the pointwise stabilizer of $\Delta.$ This contradicts
the fact that $vt \ne v.$ We conclude that the Sylow $p$-subgroups of $G$
are elementary abelian, proving (4).

Suppose, to get a contradiction, that $V$ has an imprimitivity decomposition
$V=W_1 \oplus \ldots \oplus W_k.$ We may suppose that $H$ permutes the
$W_i$ primitively. Let $M=N_H(W_1)$, so that $M$ is a maximal subgroup of
$H$. If $M$ contains $O_{p^{\prime}}(H)$, then $|H:M|=p$ and $M$ is normal in $H$.
If not, then $MO_{p^{\prime}}(H)=H$ and so $p^2$ divides $|M|$. The first
possibility is easy to handle. We have $k=p.$ We let $\Delta=W_1.$ Then 
$\mathrm{Stab}_H(\Delta)=M$ and $\mathrm{Stab}_G(\Delta)=W_1M.$
Hence $|\mathrm{Stab}_G(\Delta)|_p=|M|_p=p,$ a contradiction.

Suppose then that $|M|_p=p^2.$ Let $L=\mathrm{core}_H(M).$ Then
$M/L$ has a normal complement $K/L$ in $H/L.$ Since $O^{p{\prime}}
(H)=H,$ we see that $K/L$ is a noncentral chief factor of $H$. The 
permutation group $H/L$ on the $W_i$ is the affine group $(K/L)(M/L).$
We have $k=|K/L|.$ Suppose first that $p^2$ divides $|M/L|$. Since
$|H/L|<|G|$ and $G$ is a minimal counterexample to Theorem A, we see
that $H/L$ is a $p$-moderate primitive affine group. Hence there is a 
subset $\Delta$ of $\{W_1, \ldots, W_k\}$ whose setwise stabilizer
$S/L$ in $H/L$ satisfies $|S/L|_p=p.$ Let $W$ be the direct sum of the
$W_i$ in $\Delta$. Then $\mathrm{Stab}_G(W)=W\mathrm{Stab}_H
(\Delta).$ Since $(p,|L|)=1,$ we have
$$|\mathrm{Stab}_H(\Delta)|_p=|\mathrm{Stab}_{H/L}(\Delta)|_p=
|S/L|_p=p.$$ Thus $G$ is $p$-moderate, a contradiction.

The other possibility is that $|M/L|_p=|L|_p=p.$ We have two sub-cases.
If $H/L$ is not a $p$-concealed permutation group on $\{W_1, \ldots ,
W_k\},$ choose $\Delta \subseteq \{W_1, \ldots , W_k\}$ such that $\mathrm
{Stab}_{H/L}(\Delta)$ is a $p^{\prime}$-group. As above, let $W$ be the 
direct sum of the $W_i$ in $\Delta.$ As above $\mathrm{Stab}_G(W)=
W\mathrm{Stab}_H(\Delta).$ Then
$$|\mathrm{Stab}_H(\Delta)|_p=|\mathrm{Stab}_L(\Delta)|_p=
|L|_p=p.$$ Thus $|\mathrm{Stab}_G(W)|_p=p$ and so $G$ is 
$p$-moderate, a contradiction.

In the other sub-case, $H/L$ is $p$-concealed. Then $(p, H/L, k)$
is either $(2,D_6,3),$ or $(2,D_{10},5)$, or $(3,J,8).$ Suppose first that
$p=2.$ Let $\Delta \subseteq W_1 \oplus W_2$ be the set of all vectors in
$W_1 \oplus W_2$ whose $W_2$-component is nonzero. Now $D_{2k}$ 
is a Frobenius group on $\{W_1, \ldots ,W_k\}$ and so only the identity in
$H/L$ fixes both $W_1$ and $W_2.$ We have
$$L<\mathrm{Stab}_G(\Delta)=W_1L.$$ Hence $|\mathrm{Stab}_G(\Delta)|_2
=2$ and so $G$ is $2$-moderate, a contradiction.

Suppose then that $p=3.$ Let $\Delta \subseteq W_1 \oplus W_2 \oplus W_3 \oplus
W_4$ be the set of vectors whose $W_3$- and $W_4$- components are 
nonzero. Note that the $|W_i|$ are not equal to 2; that would imply that
$L=1$ and so $O_2(H)>1,$ against the irreducibility of $V$ as an
$H$-module. It follows that if $vh \in \mathrm{Stab}_G(\Delta)$ with
$v \in V$ and $h \in H$, then $v \in W_1 \oplus W_2$ and $h$ permutes $W_1$
and $W_2$ and $h$ also permutes $W_3$ and $W_4$. Hence if $g$, and
therefore $h$, have order 3, we must have $h \in L.$ Thus
$$|\mathrm{Stab}_G(\Delta)|_3=|\mathrm{Stab}_{VL}(\Delta)|_3=
|(W_1 \oplus W_2)L|_3=3.$$
Hence $G$ is 3-moderate, a contradiction. This proves (5). \qedsymbol

\bigskip

\begin{center}

\textbf{2. Proof of Theorem A}

\end{center}

\bigskip

\noindent \textbf{Theorem 2.1} Let $G=(VH, V)$ be a solvable affine primitive
permutation group. Suppose that $G$ satisfies conditions (1)-(5) of Theorem
1.3. Then $H$ has a regular orbit on $V \oplus V.$

\medskip

\noindent \textit{Proof.} This is an immediate consequence of Dolfi's theorem 
[2, Theorem 3.4]. Indeed, since $H$ is a primitive linear group on $V$,
[2, Theorem 3.4] says that $H$ has a regular orbit on $V \oplus V,$ apart from
four exceptional cases. In the first exceptional case, $H \cong S_3$ and $|V|=4.$ 
In the second, $H$ is $SL(2,3)$ or $GL(2,3)$ and $|V|=9.$ In the third, $H$ is 
$3^{1+2}.SL(2,3)$ or $3^{1+2}.GL(2,3)$ and $|V|=64.$ In the fourth, $|V|=81$
and $$(Q_8 \ast Q_8)(Z_3 \times Z_3) \leq H \leq (Q_8 \ast Q_8)(S_3 \times S_3).$$
In each of these cases, one checks that conditions (1)-(4) of Theorem 1.3 are
not all satisfied. \qedsymbol

\medskip

\noindent \textit{Proof of Theorem A for p greater than 2}. We will show that a minimal
counterexample to Theorem A, as described in Theorem 1.3, is in fact 
$p$-moderate, the desired contradiction.

By Theorem 2.1, some $(v,w) \in V\oplus V$ lies in a regular $H$-orbit. We  may
assume that $v$ and $w$ are both nonzero. Let $t \in H$ be a fixed element of order 
$p$. Let $\mathcal O_1=\{v,vt , \ldots, vt^{p-1}\}$ and let $\mathcal O_2=
\{w, wt, \ldots, wt^{p-1}\}$. We cannot have $|\mathcal O_1|=|\mathcal O_2|=1.$
If, say, $|\mathcal O_1|=p$ and $|\mathcal O_2|=1,$ replace $w$ by $v+w$. Then
$(v,v+w)$ lies in a regular $H$-orbit and $t$ centralizes neither $v$ nor $v+w$, so 
we may assume that $|\mathcal O_1|=|\mathcal O_2|=p.$ If $\mathcal O_1=\mathcal
O_2,$ let $\Delta=\{0\} \cup \mathcal O_1.$ Let $t \in P \in \mathrm{Syl}_p(\mathrm
{Stab}_G(\Delta)).$ Since $P$ is abelian, $P$ stabilizes $\{0\}$ and $\mathcal O_1.$
Thus $P \leq H.$ Since $(v,w)$ lies in a regular $H$-orbit, $P$ acts faithfully on
$\mathcal O_1.$ Hence $|P|=p$ and so $G$ is $p$-moderate, a contradiction.

Hence $|\mathcal O_1|=|\mathcal O_2|=p$  and $\mathcal O_1$ and $\mathcal O_2$
are disjoint. Let $\Delta =\mathcal O_1 \cup \mathcal O_2 \cup \{0\}.$ Again let
$t \in P \in \mathrm{Syl}_p(\mathrm{Stab}_G(\Delta)).$ Since $P$ is abelian and
$p \ne 2,$ $P$ must stabilize $\{0\}, \mathcal O_1,$ and $\mathcal O_2.$ Suppose 
that $|P|=p^2.$ Then $P$ acts intransitively on $\mathcal O_1 \cup \mathcal O_2$ and 
there exists $t^{\prime} \in P$ such that $t^{\prime}$ fixes every point in 
$\mathcal O_2$ and $t^{\prime}$ transitively permutes the points in $\mathcal O_1.$
Now let $\Gamma=\{0,w\} \cup \mathcal O_1$. Then $t^{\prime} \in \mathrm{Stab}_G(\Gamma)$.
Let $t^{\prime} \in Q \in \mathrm{Syl}_p(\mathrm{Stab}_G(\Gamma)).$ Then $Q$ fixes
0 and $w$, and $Q$ acts transitively on $\mathcal O_1.$ In particular, $Q \leq H.$
The definition of $(v,w)$ now implies that $|Q|=p,$ and so $G$ is $p$-moderate, a
contradiction. Hence $|P|=p,$ which also implies that $G$ is $p$-moderate, the desired
contradiction. This proves Theorem A for $p>2.$  \qedsymbol

\medskip

\noindent \textbf{Lemma 2.2} Suppose $(G,V)$ is a primitive affine group which
satisfies conditions (1)-(5) of Theorem 1.3 for $p=2.$ Suppose that $H$ has a
regular orbit on $V$. Then $(G,V)$ is 2-moderate.

\medskip

\noindent \textit{Proof.} Let $v \in V$ lie in a regular $H$-orbit. Let $t\in H$ be a 
fixed involution and let $\Gamma=\{0,v,vt\}$, so that $t\in \mathrm{Stab}_H(\Gamma).$
Suppose that $g \in G$ is an involution that commutes with $t$. If $g$ stabilizes
$\Gamma$, then $g$ must stabilize the two $\langle t \rangle$-orbits on $\Gamma.$
Thus $g \in H$ and $g$ stabilizes $\{v, vt\}$. Since $v$ lies in a regular $H$-orbit,
$g$ must interchange $v$ and $vt.$ But then $gt$ fixes $\{v, vt\}$ pointwise. It
follows that $gt=1$, and so $g=t.$ Thus $|\mathrm{Stab}_G(\Gamma)|_2=2,$
as desired. \qedsymbol

\medskip

If $H$ is not metacyclic, our goal is to show that $H$ has a regular orbit on
$V$. We accomplish this by analyzing extraspecial subgroups of $F(H)$
and applying the results of [7].

\medskip

\noindent \textbf{Proposition 2.3} Suppose $(G,V)$ is a primitive affine group
which satisfies conditions (1)-(5) of Theorem 1.3 for $p=2.$ If $H$ is not
metacyclic, then $H$ has a well-defined normal subgroup $E$, which is a
direct product of extraspecial groups, with one extraspecial direct factor
for each prime divisor of $|E|$. Then $e:=\sqrt{|E/Z(E)|}$ is odd. If $e>3,$
then $H$ has a regular orbit on $V$.

\medskip

\noindent \textit{Proof.} We follow the notation of [12, Theorem 2.1]. We
note that $E$ as above is not quite the same as the subgroup $E$
defined in the standard reference [9, Corollary 1.10]. Indeed $E$ as
above corresponds to the subgroup $D$ defined in [9, Lemma 2.10].
Since $E/Z(E)$ is $H$-isomorphic to $D/Z(D),$ the distinction between
the two subgroups is not important.

Clearly $e$ is odd; otherwise $H$ would contain an extraspecial 
2-group, contradicting $|H|_2=4.$ For $e>3,$ we consult [7, Table 4.1]. We note
that $p$ in that table denotes the prime divisor of $|V|.$ We see from Table 4.1 
that regular $H$-orbits must exist whenever $|V|$ is odd. \qedsymbol

\medskip

\noindent \textbf{Lemma 2.4} Suppose $(G,V)$ is a primitive affine group
which satisfies conditions (1)-(5) of Theorem 1.3 for $p=2$. Then $e \ne 3.$

\medskip

\noindent \textit{Proof.} We follow [9, Corollary 1.10]. We let $F$ denote the
Fitting subgroup of $H$ and we let $Z$ denote the socle of the cyclic group
$Z(F).$ We let $A=C_H(Z)$. Suppose $e=3.$ Then $F$ contains a normal
subgroup $D$ of $H$ which is extraspecial of order 27 and exponent 3, and 
$F=DC_F(D).$ By [9, Corollary 1.10(viii)], $D/Z(D)$ is a completely reducible
$H/F$-module and a faithful $A/F$-module. Hence $D/Z(D)$ is also a 
completely reducible $A/F$-module.

By [9, Corollary 1.10(ix)], $A/F$ is isomorphic to a subgroup of $SL(2,3).$
Since $|H|_2=4,$ it follows that $|A/F|$ divides 12. From the structure of
$SL(2,3)$ we see that $O_3(A/F)>1$ whenever 3 divides $|A/F|.$ Since
$D/Z(D)$ is a completely reducible $A/F$-module, we conclude that $A/F$ 
is an elementary abelian 2-group. Hence $|A/F| \leq 2.$ Since $O^{2^{\prime}}
(H)=H$ by Theorem 1.3(2) and $H/A$ is abelian, we see that $H/A$ and therefore
$H/F$ are elementary abelian 2-groups. By [9, Corollary 1.10(iii),(iv)], however,
the chief factors of $H$ within $D/Z(D)$ have square order. This contradicts
the simultaneous diagonalizability of $H/F$ on $D/Z(D)$, which proves
that $e \ne 3,$ as desired. \qedsymbol

\medskip

\noindent \textit{Proof of Theorem A for p=2}. By Lemma 2.2, Proposition 2.3,
and Lemma 2.4, we may assume that $H$ is metacyclic and $G=VH$ satisfies
conditions (1)-(5) of Theorem 1.3. Let $C$ be a normal cyclic subgroup of
$H$ such that $H/C$ is cyclic. Since $H$ has an elementary abelian Sylow
2-subgroup of order 4, it follows that $|C|_2=|H/C|_2=2.$ Hence $|O_2(H)|=2.$
Let $O_2(H)=\langle z \rangle.$

Let $u \in H$ be a noncentral involution. Let $v \in V$ be a $-1$-eigenvector for
$u$ and let $w \in V$ be a +1-eigenvector for $u$. Let $\Gamma=\{0,w,v,-v\}.$
Then $u$ stabilizes $\Gamma.$ For $y\in V$ and $h \in H,$ the permutation $yh$
in $G$ sends $a \in V$ to $(a+y)h \in V.$ Suppose that the permutation $xs \in G$ 
is an involution that commutes with $u$ and stabilizes $\Gamma.$ Then $xs$
permutes $\{0\}, \{w\},$ and $\{v, -v\}$, the orbits of $u$ on $\Gamma$. Thus $xs$
stabilizes $\{v, -v\}$ and $\{0,w\}.$

Suppose, to get a contradiction, that $xs$ interchanges 0 and $w$. Then the vector
$xs$ equals $w$, and $(w+x)s=0,$ and so $x=-w.$ Since $s$ must be an involution, we 
have $x=ws.$ Hence $ws=-w.$

Now the permutation $xs$ either interchanges $v$ and $-v$, or $xs$ fixes both $v$ and $-v.$
If the first possibility holds, then $(v-w)s=-v$ and $(-v-w)s=v$. Adding equations gives
$-2ws=0,$ which implies that $w=0,$ a contradiction. If the second possibility holds,
adding equations again leads to $w=0,$ which is false.

Therefore $xs$ fixes both 0 and $w$. Then $xs=s \in H,$ and $ws=w$. If $s \ne u,$
then $\langle u,s \rangle$ contains $O_2(H)=\langle z \rangle.$ But this is a 
contradiction, since $u$ and $s$ stabilize $\Gamma$ while $z$ does not. Therefore
$s=u$ and so $|\mathrm{Stab}_G(\Gamma)|_2=2.$ Thus $(G,V)$ is 2-moderate,
the final contradiction. \qedsymbol

\medskip

\noindent \textit{Remark.} The referee has pointed out that Theorem A becomes false
if one substitutes ``transitive" for ``primitive" in its statement. An example is
$\mathrm{Sym}(3) \wr \mathrm{Alt}(3)$ acting imprimitively on nine points, for $p=2.$

\bigskip

\begin{center}

\textbf{3. Non-solvable groups}

\end{center}

\bigskip

We saw in the preceding sections that the problem of classifying the $p$-moderate
solvable primitive permutation groups $(G, \Omega)$ is much more challenging
when $G$ has elementary abelian Sylow $p$-subgroups. The next result suggests 
that the same is true when $G$ is not solvable.

\medskip

\noindent \textbf{Proposition 3.1} Let $(G,\Omega)$ be a (not necessarily primitive)
permutation group. Let $P \in \mathrm{Syl}_p(G)$ and suppose that $P$ is not
elementary abelian. Let $z$ be an element of order $p$ in $\Phi(P) \cap Z(P).$
Let $f$ be the number of fixed points of $z$ on $\Omega,$ and let $n=|\Omega|$. If
$$|G:N_G(P)|<2^{(n-f)(1/p-1/p^2)},$$ then $(G, \Omega)$ is $p$-moderate.

\medskip

\noindent \textit{Proof.} Suppose $\mathcal O$ is a $P$-orbit in $\Omega$
with $|\mathcal O|=m<p^2$. Since a Sylow  $p$-subgroup of the symmetric group
$S_m$ is elementary abelian, $z$ must act trivially on $\mathcal O.$ Hence if we 
let $\Omega^{\prime}$ denote the set of points in $\Omega$ that are not fixed 
by $z$, then each $P$-orbit on $\Omega^{\prime}$ has size at least $p^2.$
It follows that $P$  has at most $f+(n-f)/p^2$ orbits on $\Omega,$ while $z$ has exactly
$f+(n-f)/p$ orbits on $\Omega.$ Hence $z$ stabilizes exactly $2^{f+(n-f)/p}$ subsets of 
$\Omega$ and the number of subsets of $\Omega$ stabilized by some Sylow
$p$-subgroup of $G$ is at most $|G:N_G(P)|2^{f+(n-f)/p^2}.$ Hence $G$ is 
$p$-moderate if
$$|G:N_G(P)|2^{f+(n-f)/p^2}<2^{f+(n-f)/p}.$$
This inequality is equivalent to the one in the statement of the proposition. \qedsymbol

\newpage

\begin{center}
\textbf{References}
\end{center}
\bigskip
\noindent  \textbf1. L. Babai, Asymmetric coloring of locally finite graphs
and profinite permutation groups, J. Algebra 607(2022), 64-106.

\medskip

\noindent \textbf2. S. Dolfi, Large orbits in coprime actions of solvable groups,
Trans. Amer. Math. Soc. 360(2008), 135-152.

\medskip

\noindent \textbf3. C. Eaton and A. Moret\'o, Extending Brauer's height zero conjecture
to blocks with nonabelian defect groups, Int. Math. Res. Not. IMRN 20(2014), 5581-5601.

\medskip

\noindent \textbf4. M. Giudici, M. Liebeck, C. Praeger, J. Saxl, and P. H. Tiep,
Arithmetic results on orbits of linear groups, Trans. Amer. Math. Soc.
368(2016), 2415-2467.

\medskip

\noindent \textbf5. D. Gluck, Set-stabilizers in solvable permutation groups, J.
Group Theory 28(2025), 1297-1305.

\medskip

\noindent \textbf6. D. Gluck and T. R. Wolf, Brauer's height conjecture for
$p$-solvable groups, Trans. Amer. Math. Soc. 282(1984), 137-152.

\medskip

\noindent \textbf7. D. Holt and Y. Yang, Regular orbits of finite primitive
solvable groups, the final classification, J. Symbolic Comput. 116(2023), 139-145.

\medskip

\noindent \textbf8. G. Malle, G. Navarro, A. Schaeffer Fry, P. H. Tiep,
Brauer's height zero conjecture, Ann. of Math. 200(2024), 557-608.

\medskip

\noindent \textbf9. O. Manz and T. R. Wolf, Representations of Solvable Groups,
Cambridge University Press, Cambridge, 1993.

\medskip

\noindent \textbf{10}. G. Navarro and P. H. Tiep, Characters of relative 
$p^{\prime}$-degree over normal subgroups, Ann. of Math. 178(2013),
1135-1171.

\medskip

\noindent \textbf{11}. L. Sabatini, On stabilizers in finite permutation groups, Bull. London Math. Soc. 58(1), 
2026.

\medskip

\noindent \textbf {12}. A. Vasil'ev, E. Vdovin, Y. Yang, Regular orbits of finite
primitive solvable groups, III, J. Algebra 590(2022), 139-154.

\end{document}